\documentclass[10pt]{amsart}
\usepackage{amsmath,amssymb,amsthm}
\usepackage{combelow}
\usepackage[colorlinks,hyperindex,linkcolor=blue]{hyperref}
\usepackage[utf8]{inputenc}
\usepackage[margin=1.2in]{geometry}

\newcommand{\diffto}{\xrightarrow{\raisebox{-0.2 em}[0pt][0pt]{\smash{\ensuremath{\sim}}}}}

\newtheorem*{example}{Example}
\newtheorem*{remark}{Remark}
\newtheorem*{definition}{Definition}
\newtheorem{proposition}{Proposition}
\newtheorem{theorem}[proposition]{Theorem}
\newtheorem*{corollary}{Corollary}
\newtheorem*{mtheorem}{Main Theorem}
\newtheorem*{stheorem}{Theorem}
\newtheorem{question}{Question}

\newcommand{\Addresses}{{
  \bigskip
  \footnotesize

 Florian Zeiser,\par\vspace*{0pt}\nopagebreak
\textsc{University of Illinois at Urbana-Champaign, 61801 Urbana, United States}\par\vspace*{-2pt}\nopagebreak
\textit{E-mail address}: \texttt{fzeiser@illinois.edu}
}}

\title{A remark on the linearization of nondegenerate Nambu structures of coorder 1}
\author{Florian Zeiser}

\begin{document}

\maketitle
\begin{abstract}
    In this note we show that Nambu structures of coorder 1 can always be linearized if they admit a closed integrable differential form. In particular, we show that a unimodular Poisson structure whose isotropy Lie algebra at a singular point is $\mathfrak{sl}(2,\mathbb{R})$, can always be linearized. 
\end{abstract}
Nambu structures were introduced by Nambu \cite{Nambu} and expanded on by Takhtajan \cite{Tak}, in an attempt to generalize Hamiltonian mechanics. A Nambu structure of order $q$ on a smooth manifold $M$ is given by a $q$-vector field $\Pi\in \mathfrak{X}^q(M)$ which is preserved by all its Hamiltonian vector fields, i.e. 
\[ \mathcal{L}_{X_{f_1,\dots ,f_{q-1}}}\Pi=0 \qquad \text{ where } \ X_{f_1,\dots ,f_{q-1}}:= \Pi(\mathrm{d}f_1, \dots ,\mathrm{d}f_{q-1},\cdot) \ \text{ and } f_1,\dots,f_{q-1}\in C^\infty (M).\]
We call $p:=\dim(M)-q$ the coorder of $\Pi$. In particular, we recover the definition of a Poisson structure for $q=2$. Hence, they are also known as Nambu-Poisson structures. To describe normal forms of Nambu structures we distinguish the cases $q=2$ and $q\ge 3$. The case of vector fields, i.e. when $q=1$, is not treated here. For more details on Nambu structures we refer to \cite{DZ}.\vspace*{5pt}


\underline{\textbf{The Poisson case:}}\vspace*{2pt}

Weinstein showed in \cite{Wein83} that, locally around a point, any Poisson structure decomposes as product of a symplectic manifold and a manifold with vanishing Poisson structure.  Hence the existence of local normal forms for Poisson structures reduces to the question of linearization, i.e. given a Poisson structure $\pi$ on $M$ and $p\in M$ with $\pi(p)=0$, the isotropy Lie algebra $\mathfrak{g}_p$ on $T_p^*M$ is defined by
\[ [\mathrm{d}_{p}g,\mathrm{d}_{p}h ]:=\mathrm{d}_{p}\pi(\mathrm{d} g,\mathrm{d} h),\quad \quad g,h\in C^{\infty}(M).\] 
Any Lie algebra $(\mathfrak{g} ,[\cdot,\cdot])$ admits an associated Poisson structure $\pi_{\mathfrak{g}}$ on the dual space $\mathfrak{g}^*$ defined by
\begin{align*}
      \pi_{\mathfrak{g},\xi}(X,Y):= \xi([X,Y]) \quad \text{where} \ \xi \in \mathfrak{g}^*, \ X,Y\in \mathfrak{g}. 
\end{align*}
Applying this construction to the isotropy Lie algebra $\mathfrak{g}_p$ of $\pi$ at $p\in M$, one obtains the linear Poisson structure $\pi_{\mathfrak{g}_{p}}$ on $\mathfrak{g}^*_{p}=T_{p}M$, which plays the role of the first order approximation of $\pi$ at $p$.

\begin{definition}
A Poisson structure $\pi$ on $M$ is called \textbf{linearizable} around a zero $p\in M$ of $\pi$ if there exists a smooth Poisson diffeomorphism
\[ \Phi: (O,\pi)\diffto (V,\pi_{\mathfrak{g}_{p}}),\quad \textrm{with}\quad \Phi(p)=0,\]
for some open neighborhoods $p\in O\subset M$ and $0\in V\subset  \mathfrak{g}^*_{p}$.
\end{definition}

The question of linearization has been of interest ever since Weinstein's work (see \cite{Conn85,Daz,Wein87,MZ04,DZ,CF11,,MZ221}). Of particular interest for this work is the following:
\begin{proposition}[Weinstein \cite{Wein83}]\label{th:sl2r}
    There exists a Poisson structure $\pi$ on $\mathbb{R}^3$ with $\pi(0)=0$ and $\mathfrak{g}_0=\mathfrak{sl}_2(\mathbb{R})$ which is not linearizable.
\end{proposition}
We refer to \cite{MZ221} for a more detailed outline and references on the linearization problem.\vspace*{5pt}

\underline{\textbf{The Nambu case:}}\vspace*{2pt}

If the order $q$ of the Nambu structure $\Pi$ is at least three, then $\Pi$ induces a singular foliation by leaves of dimension $q$ and $0$ respectively, and $\Pi$ is decomposable around regular points (Alekseevsky \& Guha \cite{AG}, Gautheron \cite{Gaut}, Nakanishi \cite{Nakanishi_98}).

Around singular points, i.e. $p\in M$ with $\Pi(p)=0$, we can ask again the question of linearization. In this case $\Pi$ gives rise to an $q$-ary Lie algebra $\mathfrak{F}$ (or Filippov algebra \cite{Fil}) on $T^*_pM$ by
\[ [\mathrm{d}_pf_1,\dots, \mathrm{d}_pf_q] := \mathrm{d}_p \Pi (\mathrm{d}f_1, \dots, \mathrm{d}f_q), \qquad f_1,\dots f_q\in C^\infty(M).\]
Contrary to the Lie algebra case, not every $q$-ary Lie algebra (for $q\ge 3$) gives rise to a Nambu structure on its dual. However, those induced by Nambu structures do and we obtain a linear Nambu structure $\Pi_l$ on $T_pM$ by
\begin{align*}
      \Pi_{l}(\alpha_1,\dots, \alpha_q)(Z):= [\alpha_1,\dots, \alpha_q](Z) \quad \text{where} \ Z \in T_pM, \ \alpha_1,\dots \alpha_q\in T_p^*M. 
\end{align*}
As in the Poisson case we can ask, for a given Nambu structure $\Pi$ and $p\in M$ with $\Pi(p)=0$, does there exist a smooth diffeomorphism of Nambu structures 
\[ \Phi : (O,\Pi )\to (V,\Pi_{l})\]
between open neighborhoods of $p\in M$ and $0\in T_pM$, respectively?
\begin{definition}
    If such a $\Phi$ exists, then we call $\Pi$ linearizable at $p\in M$. 
\end{definition}
Linear Nambu structures have been classified by Dufour \& Zung \cite{DZ99} and Medeiros \cite{Med}. Both make use of a dual notion for Nambu structures due to Dufour \& Zung.
\begin{proposition}[Dufour \& Zung \cite{DZ99}]\label{prop:dual}
    Let $\Pi\in \mathfrak{X}^q(M)$ and $\mu\in \Omega^m(M)$ a volume form on $M$. Set
    \[ \omega:=\iota_\Pi \mu.\]
    Then $\Pi$ is Nambu if and only if $\omega$ satisfies for all $\xi\in \mathfrak{X}^{m-q-1}(M)$ the equations
    \begin{equation}\label{eq:coN}
        \iota_\xi \omega \wedge \omega=0 \qquad \text{ and }\qquad \iota_\xi\omega\wedge \mathrm{d} \omega=0.
    \end{equation}
\end{proposition}
\begin{definition}
    We call $\omega\in\Omega^p(M)$ satisfying \eqref{eq:coN} an integrable differential form of order $p$.
\end{definition}
Note that this is coherent with the notion of coorder, i.e. the coorder of a Nambu structure is the order of an associated integrable differential form.
\begin{remark}
Integrable differential forms have been studied in the context of first integrals and singular foliations (see \cite{WR52, Ku64, Mou762, Mou76, Me76, Ca78, Mou83, Med}).     
\end{remark}
By Dufour \& Zung \cite{DZ99} and Medeiros \cite{Med}, every linear Nambu structure $\Pi_l$ of order $q\ge 3$ on an $m$-dimensional vector space $V$ belongs to one of the following two types
    \begin{itemize}
        \item \textbf{Type 1:}
        \[ \Pi_l = \sum_{i=1}^{r+1}\pm x_i \partial_{x_1}\wedge \dots \wedge \widehat{\partial_{x_i}}\wedge \dots \wedge \partial_{x_{q+1}} + \sum_{i=1}^s \pm x_{q+1+i}\partial_{x_1}\wedge \dots \wedge \widehat{\partial_{x_{r+i}}}\wedge \dots \wedge \partial_{x_{q+1}} \]
        with $0\le r \le q+1$ and $0\le s\le \min(m-q-1,q+1-r)$;
        \item \textbf{Type 2:}
        \[ \Pi_l = \partial_{x_1}\wedge \dots \wedge \partial_{x_{q-1}}\wedge \Big(\sum_{i,j=q}^m b_i^jx_j\partial_{x_i}\Big)\qquad \text{ for }\quad b_i^j\in \mathbb{R}.\]
    \end{itemize}
In this article, we focus on linear Nambu structures of Type 1. We call a linear Nambu structure of Type 1 nondegenerate if $r=q+1$. Using the standard volume form $\mu_{std}= \mathrm{d}x_1\wedge \dots \wedge \mathrm{d}x_m$ on $\mathbb{R}^m$, we obtain for the associated integrable differential form the expression
\begin{equation}\label{eq:nondeg} \iota_{\Pi_l}\mu_{std}= \mathrm{d}x_{q+2}\wedge \dots \wedge \mathrm{d}x_m \wedge \mathrm{d} f \qquad \text{ where } \quad f= \frac{1}{2}\Big(\sum_{i=1}^l x_i^2-\sum_{i=l+1}^{q+1}x_i^2\Big)
\end{equation}
up to permutation. The signature $(l,q+1-l)$ of $f$ is, up to permutation, a discrete invariant of $\Pi$.

\begin{theorem}[Zung \cite{Zung}]\label{thm:linearN}
    Let $\Pi$ be a Nambu structure on $M$ of order $q\ge 3$ and $p\in M$ such that $\Pi(p)=0$. If the linear part is nondegenerate of Type 1 with signature different from $(2,\star)$ and $(\star,2)$, then $\Pi$ is linearizable at $p$.
\end{theorem}
For Nambu structures whose linear part is nondegenerate of Type 1 of signature $(2,\star)$ and $(\star,2)$ Dufour \& Zung \cite[Section 5]{DZ99} gave examples of coorder 2 which are not linearizable. For the coorder 1 case there are even earlier counterexamples for linearizablity, using Proposition \ref{prop:dual}.
\begin{example}[\cite{Mou762}]
    Consider $\alpha \in \Omega^1(\mathbb{R}^n)$ given by
    \[ \alpha := \mathrm{d} f +\frac{g(f)}{x_1^2+x_2^2}(x_2\mathrm{d}x_1-x_1\mathrm{d}x_2)\]
    where $f$ is the function from $\eqref{eq:nondeg}$ with $l=2$ and $g\in C^\infty(\mathbb{R})$ with $g(x)=0$ for $x\le 0$ and $g(x)>0$ else. Note that $\alpha$ satisfies \eqref{eq:coN} and hence induces a Nambu structure $\Pi$ of coorder 1 with
    \[ X=X_{x_3,\dots, x_n}= \pm \Big( x_1\frac{\partial}{\partial_{x_2}}-x_1\frac{\partial}{\partial_{x_2}}+ \frac{l(f)}{x_1^2+x_2^2}\big(x_1\frac{\partial}{\partial_{x_1}}+x_2\frac{\partial}{\partial_{x_2}}\big)\Big)\]
    The Hamiltonian vector field $X$ admits flow lines spiraling towards the cone $C:=\{f=0\}$, hence inducing holonomy at the leaf $C\setminus \{0\}$ of the regular foliation on $\mathbb{R}^n\setminus \{0\}$. Since the foliation induced by the linear Nambu structure does not have holonomy, there can not exist a linearization. For $n=3$, under a duality as in Proposition \ref{prop:dual}, this corresponds to Weinstein's example to show Proposition \ref{th:sl2r}.
\end{example}
Here we show that Nambu structures which admit a closed integrable 1-form are always linearizable.
\begin{mtheorem}
    Let $\Pi$ be a Nambu structure on $M^m$ of order $m-1$ with $m\ge 3$ and $p\in M$ such that $\Pi(p)=0$. Assume the linearization $\Pi_{l}$ is nondegenerate of type 1 and there exists a volume form $\mu$ around $p\in M$ such that
    \begin{equation}\label{eq:unimodularN}
        \mathrm{d}\iota_\Pi \mu=0.
    \end{equation}
    Then $\Pi$ is linearizable at $p$.
\end{mtheorem}
The Main theorem, along with Theorem \ref{thm:linearN}, leads to the question:
\begin{question}
    Let $\Pi$ be a Nambu structure on $M$ of order $q\ge 3$ and $p\in M$ such that $\Pi(p)=0$. Assume the linear part is nondegenerate of Type 1 with signature $(2,\star)$ or $(\star,2)$ and there exists a volume form $\mu$ around $p\in M$ such that
    \begin{equation*}
        \mathrm{d}\iota_\Pi \mu=0.
    \end{equation*}
    Is $\Pi$ always linearizable at $p$ in this case?
\end{question}
For the Poisson case, the condition \eqref{eq:unimodularN} is well-known. Recall from \cite{Wein97} that a Poisson structure $\pi$ on an orientable manifold $M^m$ is called unimodular if there exists a volume form $\mu \in \Omega^m(M)$ with
\begin{equation*}
    \mathrm{d}\iota_\pi\mu=0.
\end{equation*}
Hence for $m=3$ the Main Theorem above yields:
\begin{corollary}
    Let $\pi$ be Poisson structure on $M^3$ which is unimodular locally around a zero $p\in M$ with $\mathfrak{g}_p\simeq \mathfrak{sl}_2(\mathbb{R})$. Then $\pi$ is linearizable around $p$.
\end{corollary}
Taking into account the isomorphism $\mathfrak{sl}_2(\mathbb{R})\simeq \mathfrak{su}(1,1)$ and the examples of non-linearizable Poisson structures for semisimple Lie algebras of real rank 1 with non-semisimple compact part by Monnier \& Zung \cite{DZ99}[Theorem 4.3.9] one is inclined to ask:
\begin{question}
    Let $\pi$ be a Poisson structure on $M$ and $p\in M$ such that $\pi(p)=0$. Assume the isotropy Lie algebra $\mathfrak{g}_p$ is semisimple of real rank 1 with non-semisimple compact part. If $\pi$ is unimodular around $p$, is $\pi$ linearizable at $p$?
\end{question}

The proof of the Main Theorem is obtained in two steps. We first linearize the singular foliation using an isochore Morse Lemma due to Colin de Verdi\'ere \& Vey \cite{VV79}. Then we use a Moser argument to linearize the associated Nambu structure. We want to point out that the linearization of the associated singular foliation was already shown by Moussu \cite[Theorem C]{Mou762} under the condition that the foliation is without holonomy.

\subsection*{Acknowledgements} I would like to thank Ioan Marcut for making me aware of the reference \cite{VV79}.

\subsection*{Proof of the Main Theorem}
   Let $\mu \in \Omega^n(\mathbb{R}^n)$ be a volume form such that 
    \[ \mathrm{d}\iota_\Pi \mu= 0.\]
    Since $\mathbb{R}^n$ has trivial first de-Rham cohomology, there exists a function $g\in C^\infty(\mathbb{R}^n)$ such that
    \[ \mathrm{d}g=\iota_\Pi \mu \qquad \text{ and }\qquad g(0)= 0.\]
    Note that $\Pi$ is uniquely determined by the pair $(g,\mu)$. By assumption, the second order Taylor polynomial of $g$ gives, up to a linear change of coordinates, rise to a non-degenerate quadratic form $f$ on $\mathbb{R}^n$ as described in \eqref{eq:nondeg} for $q+1=n$. 
    \begin{stheorem}[\cite{VV79}]
        There exist open neighborhoods $O,O'\subset \mathbb{R}^n$ and $U\subset \mathbb{R}$ of the origins respectively, a smooth function $h\in C^{\infty}(U)$ and a diffeomorphism $\phi : O\to O'$ such that
        \[ (f,h(f)\mu_{std})= \phi^*(g,\mu).\] 
    \end{stheorem}
    
    Since $\Pi$ is uniquely determined by the pair $(f,h(f)\mu_{std})$ we have 
    \[ \Pi = k(f)\Pi_{l}\]
    where $k(f)=h(f)^{-1}$, and due to the condition on the first jet we have $k(0)=1$. We finish the proof using a Moser argument. For $t\in [0,1]$ consider the family of $(n-1)$-vector fields given by
    \[ \Pi_t= t\Pi_l +(1-t)\Pi_{l}= \big(1 +t(k(f)-1)\big)\Pi_{l}.\]
    Note that this is a family of Nambu structures since $[\Pi_l,f]=0$. Denote by $E$ be the Euler vector field on $\mathbb{R}^n$ and $r\in C^{\infty}(\mathbb{R}\times[0,1])$. We make the Ansatz:
    \[ \Phi_t^*(\Pi_t)=\Pi_l, \]
    where $\Phi_t$ is the time-$t$ flow with respect to the vector field $r_t(f)E$. From the Ansatz we obtain
    \[ 0=\mathcal{L}_{r_t(f)E}\Pi_t+\partial_t\Pi_t= \biggl(r_t(f)\Bigl(t\bigl(2f\partial_x(k)(f)-(n-2)(k(f)-1)\bigr)-(n-2)\Bigr) +(k(f)-1)\biggr)\Pi_l.\]
    Sufficiently close to the origin in $\mathbb{R}^n$ the function 
    \[ r_t(f):= \frac{k(f)-1}{(n-2)(1+t(1-k(f)))-2tf\partial_x(k)(f))}\]
    is well-defined for all $t\in [0,1]$. Hence the flow $\Phi_t$ of $r_t(f)E$ is well defined for all $t\in [0,1]$ sufficiently close to the origin and satisfies
    \[ \Phi_t^*\Pi_t= \Pi_{l}.\]
    Therefore $\Phi_1$ maps $\Pi$ to $\Pi_l$.

\bibliographystyle{unsrt}
\bibliography{main}

\begin{thebibliography}{10}

\bibitem{Nambu}
Y.~Nambu.
\newblock Generalized {H}amiltonian dynamics.
\newblock {\em Phys. Rev. D (3)}, 7:2405--2412, 1973.

\bibitem{Tak}
L.~Takhtajan.
\newblock On foundation of the generalized {N}ambu mechanics.
\newblock {\em Comm. Math. Phys.}, 160(2):295--315, 1994.

\bibitem{DZ}
J.-P. Dufour and N.~T. Zung.
\newblock {\em Poisson structures and their normal forms}, volume 242 of {\em
  Progress in Mathematics}.
\newblock Birkh\"{a}user Verlag, Basel, 2005.

\bibitem{Wein83}
A.~Weinstein.
\newblock The local structure of {P}oisson manifolds.
\newblock {\em J. Differential Geom.}, 18(3):523--557, 1983.

\bibitem{Conn85}
J.~F. Conn.
\newblock Normal forms for smooth {P}oisson structures.
\newblock {\em Ann. of Math. (2)}, 121(3):565--593, 1985.

\bibitem{Daz}
P.~Dazord.
\newblock Stabilit\'{e} et lin\'{e}arisation dans les vari\'{e}t\'{e}s de
  {P}oisson.
\newblock In {\em Symplectic geometry and mechanics ({B}alaruc, 1983)}, Travaux
  en Cours, pages 59--75. Hermann, Paris, 1985.

\bibitem{Wein87}
Alan Weinstein.
\newblock Poisson geometry of the principal series and nonlinearizable
  structures.
\newblock {\em J. Differential Geom.}, 25(1):55--73, 1987.

\bibitem{MZ04}
P.~Monnier and N.~T. Zung.
\newblock Levi decomposition for smooth {P}oisson structures.
\newblock {\em J. Differential Geom.}, 68(2):347--395, 2004.

\bibitem{CF11}
M.~Crainic and R.~L. Fernandes.
\newblock A geometric approach to {C}onn's linearization theorem.
\newblock {\em Ann. of Math. (2)}, 173(2):1121--1139, 2011.

\bibitem{MZ221}
I.~M\u{a}rcu\c{t} and F.~Zeiser.
\newblock The poisson linearization problem for $\mathfrak{sl}_2(\mathbb{C})$.
  part i: Poisson cohomology, 2022.

\bibitem{AG}
D.~Alekseevsky and P.~Guha.
\newblock On decomposability of {N}ambu-{P}oisson tensor.
\newblock {\em Acta Math. Univ. Comenian. (N.S.)}, 65(1):1--9, 1996.

\bibitem{Gaut}
P.~Gautheron.
\newblock Some remarks concerning {N}ambu mechanics.
\newblock {\em Lett. Math. Phys.}, 37(1):103--116, 1996.

\bibitem{Nakanishi_98}
N.~Nakanishi.
\newblock On {N}ambu-{P}oisson manifolds.
\newblock {\em Rev. Math. Phys.}, 10(4):499--510, 1998.

\bibitem{Fil}
V.~T. Filippov.
\newblock {$n$}-{L}ie algebras.
\newblock {\em Sibirsk. Mat. Zh.}, 26(6):126--140, 191, 1985.

\bibitem{DZ99}
J.-P. Dufour and N.~T. Zung.
\newblock Linearization of {N}ambu structures.
\newblock {\em Compositio Math.}, 117(1):77--98, 1999.

\bibitem{Med}
A.~S. de~Medeiros.
\newblock Singular foliations and differential {$p$}-forms.
\newblock {\em Ann. Fac. Sci. Toulouse Math. (6)}, 9(3):451--466, 2000.

\bibitem{WR52}
W.-T. Wu and G.~Reeb.
\newblock {\em Sur les espaces fibr{\'e}s et les vari{\'e}t{\'e}s
  feuillet{\'e}es}.
\newblock Institut de math{\'e}matique (Strasbourg), 1952.

\bibitem{Ku64}
I.~Kupka.
\newblock The singularities of integrable structurally stable {P}faffian forms.
\newblock {\em Proc. Nat. Acad. Sci. U.S.A.}, 52:1431--1432, 1964.

\bibitem{Mou762}
R.~Moussu.
\newblock Existence d'int\'egrales premi\`eres pour un germe de forme de
  {P}faff {$C\sp{\infty }$}, non d\'eg\'en\'er\'e.
\newblock {\em Bol. Soc. Brasil. Mat.}, 7(2):111--120, 1976.

\bibitem{Mou76}
R.~Moussu.
\newblock Sur l'existence d'int\'egrales premi\`eres pour un germe de forme de
  {P}faff.
\newblock {\em Ann. Inst. Fourier (Grenoble)}, 26(2):xi, 171--220, 1976.

\bibitem{Me76}
A.~S. de~Medeiros.
\newblock Structural stability of integrable differential forms.
\newblock In {\em Geometry and topology ({P}roc. {III} {L}atin {A}mer. {S}chool
  of {M}ath., {I}nst. {M}at. {P}ura {A}plicada {CNP}q, {R}io de {J}aneiro,
  1976)}, volume Vol. 597 of {\em Lecture Notes in Math.}, pages 395--428.
  Springer, Berlin-New York, 1977.

\bibitem{Ca78}
C.~Camacho.
\newblock Structural stability theorems for integrable differential forms on
  {$3$}-manifolds.
\newblock {\em Topology}, 17(2):143--155, 1978.

\bibitem{Mou83}
R.~Moussu.
\newblock Classification {$C\sp{\infty }$}\ des \'equations de {P}faff
  compl\`etement int\'egrables \`a{} singularit\'e{} isol\'ee.
\newblock {\em Invent. Math.}, 73(3):419--436, 1983.

\bibitem{Zung}
N.~T. Zung.
\newblock New results on the linearization of {N}ambu structures.
\newblock {\em J. Math. Pures Appl. (9)}, 99(2):211--218, 2013.

\bibitem{Wein97}
A.~Weinstein.
\newblock The modular automorphism group of a {P}oisson manifold.
\newblock {\em J. Geom. Phys.}, 23(3-4):379--394, 1997.

\bibitem{VV79}
Y.~Colin~de Verdi\`ere and J.~Vey.
\newblock Le lemme de {M}orse isochore.
\newblock {\em Topology}, 18(4):283--293, 1979.

\end{thebibliography}
\Addresses
\end{document}